\documentclass[12pt]{article}
\usepackage{amsfonts}
\usepackage{amssymb}
\usepackage{amsmath}
\usepackage{amsthm}
\usepackage{amsfonts}
\usepackage[T1]{fontenc}
\usepackage[utf8]{inputenc} 
\usepackage{cite} 
\usepackage{stmaryrd}
\usepackage[lflt]{floatflt}
\usepackage{graphicx}
\usepackage{fancyhdr}
\usepackage{color}
\usepackage[affil-sl]{authblk}
\usepackage{cancel}
\usepackage{subfigure}
\usepackage{makeidx}
\usepackage{latexsym}
\usepackage{hyperref} 
\usepackage{epsfig}
\usepackage{stmaryrd}
\usepackage{float}
\usepackage{tabularx}
\usepackage{verbatim}
\usepackage{color}
\usepackage{verbatim}
\usepackage{enumerate}
\usepackage[usenames,dvipsnames]{xcolor}

\usepackage{tabularx}

\usepackage{geometry}

\usepackage[usenames,dvipsnames]{xcolor}

\newtheorem{nota}{Note}
\newtheorem{teo}{Theorem}
\newtheorem{defi}{Definition}
\newtheorem{propo}{Proposition}
\newtheorem{lema}{Lemma}

\newtheorem{obs}{Remark}

\newcommand{\CC}{{\mathcal C}}
\newcommand{\bbR}{{\mathbb R}}

\newcommand{\e}{\epsilon}
\newcommand{\al}{\alpha}

\newcommand{\p}{\partial}

\newcommand{\G}{\Gamma}

\newcommand{\dd}{{\rm d}}

\begin{document}
\begin{center}
\LARGE
\textbf{A New Mathematical Formulation for a  {P}hase  {C}hange {P}roblem with a Memory Flux}
\end{center}

                   \begin{center}
                  {\sc  Sabrina D. Roscani, Julieta Bollati, and Domingo A. Tarzia}\\
 CONICET - Depto. Matem\'atica,
FCE, Univ. Austral,\\
 Paraguay 1950, S2000FZF Rosario, Argentina \\
(sabrinaroscani@gmail.com, jbollati@austral.edu.ar, dtarzia@austral.edu.ar)
                   \vspace{0.2cm}

       \end{center}
      
\small

\noindent \textbf{Abstract: }

  A mathematical formulation for a one-phase change problem in a form of Stefan problem with a memory flux is obtained. The hypothesis that the integral of weighted backward  fluxes is proportional to the  gradient of the temperature is considered.  The model that arises involves fractional derivatives with respect to time both  in the sense of Caputo and  of Riemann--Liouville. An integral relation for the free boundary, which is equivalent to the ``fractional Stefan condition'', is also obtained. 
		

\noindent \textbf{AMS:} Primary: 35R35, 26A33, 35C05. Secondary: 33E20, 80A22.  \\

\noindent \textbf{Note:} This paper is now published (in revised form) in Chaos, Solitons and Fractals \textbf{116} (2018), p.p. 340-347, DOI:10.1016/j.chaos.2018.09.023, and is available online at \linebreak \url{ www.elsevier.com/locate/chaos}, so always cite it with the journal's coordinates.


\section{Introduction}

The theory related to heat diffusion has been extensively developed in the last century. Modelling  classical heat diffusion comes hand in hand with Fourier Law.  Nevertheless, we shall not forget that this famous law is an experimental phenomenological principle.

In the past 40 years, many generalized flux models of the classical one (i.e. the one derived from Fourier Law) were proposed in the literature and accepted by the scientific community. See e.g.   \cite{Cha:1998,HeIg:1999,IgOs:2010,JoPr:1989,PK:1976}.

In this paper a  phase change problem  for heat  diffusion  under the hypothesis that the heat flux is a flux with memory is analysed. This kind of problems are known in the literature as Stefan problems \cite{Tarzia, Tarzia:biblio}. 

The model obtained under the memory assumption is known as an   anomalous diffusion model, and it is governed by fractional  diffusion equations. There is a vast literature in the subject of fractional diffusion equations. We refer the reader to \cite{FM-Libro, Povstenko, Pskhu-Libro} and references therein.

The study of anomalous diffusion has its origins in the investigation of non-Brownian motions (Random walks). In that context it was observed that ``the mean square displacement'' of the particles is  proportional to a power of the  time, instead of being proportional just to time. An exhaustive work in this direction has been done by Metzler and Klafter \cite{MK:2000}. Other articles in this direction are \cite{KlSo:2005, MeGlNo:1994,Nig:1984,Pa:2013}. It is worth mentioning that  many works (see e.g. \cite{BaFr:2005,GKS,Sax:1994}) suggest that the anomalous diffusion is caused by heterogeneities in the domain.

Before presenting the problem, let us establish  some usual notation related to heat conduction with the corresponding physical dimensions. 
Let us write $\textbf{T}$ for temperature, $\textbf{t}$ for time, $\textbf{m}$ for mass and $\textbf{X}$ for position.

\begin{equation}\label{medidas}
\begin{array}{ccc}
u & \textsl{temperature} & [\textbf{T}]\\
k & \textsl{thermal conductivity} & \left[\frac{{\textbf{m X}}}{{\textbf{Tt}^3}}\right] \\
\rho & \textsl{mass density} &  \left[\frac{{\textbf {m }}}{{\textbf {X}^3}}\right]  \\
c & \textsl{specific heat} &   \left[\frac{\textbf{X}^2 }{\textbf{T}\textbf{t}^2}\right]\\
d=\frac{k}{\rho c} & \textsl{ diffusion coefficient} &  \left[\frac{\textbf{X}^2}{\textbf{t}}\right]\\
l & \textsl{latent heat per  unit  mass} & \left[\frac{\textbf{X}^2 }{\textbf{t}^2}\right]
\end{array}
\end{equation}

\noindent Consider  a temperature function  $u=u(x,t)$ and  its corresponding flux $J(x,t)$, both  defined  for a semi-infinite unidimensional material. From the First Principle of Thermodynamics, we deduce the continuity equation
\begin{equation}\label{cont eq}
\rho c \frac{\p u}{\p t}(x,t)=-\frac{\p J}{\p x}(x,t).
\end{equation}
The aim of this work is to derive a model  by considering a special non-local memory flux. For example, Gurtin and Pipkin \cite{GuPi:1968} (experts in continuum mechanics and heat transfer) proposed in  1968 a general theory of heat conduction with finite velocity waves  through the following non local flux law:

\begin{equation}\label{Gen-FourierLaw}
J(x,t)=K(t)\ast \left( -k\frac{\p}{\p x}u(x,t) \right)= -k\int_{-\infty}^t K(t-\tau)\frac{\p}{\p x}u(x,\tau)\rm{d}\tau,
\end{equation}
where $K$ is a positive decreasing kernel which verifies $K(s)\rightarrow  0$  when $s\rightarrow \infty$.

Let us comment on some different explicit and implicit definitions of fluxes, and their effects on the resulting governing equations:
\begin{itemize}
\item  \textbf{Explicit forms for the flux: $J(x,t)=F(x,t)$}\\

The classical law for the flux is the \textit{Fourier Law}, which states that the flux J is proportional to the temperature gradient, that is:
\begin{equation}\label{FourierLaw}
J(x,t)=-k\frac{\p}{\p x}u(x,t).
\end{equation}

If alternatively  suppose that the flux at the point $(x,t)$ is proportional to the total flux, then the given law  is the following
\begin{equation}\label{Total-Flux}
J(x,t)=\frac{1}{\tilde{\tau}}\int_{-\infty}^t -k\frac{\p}{\p x}u(x,\tau) \dd\tau.
\end{equation}

In  (\ref{Total-Flux}), $\tilde{\tau}$ is a constant whose physical dimension is time. Another interesting thing is that (\ref{Total-Flux})  can be interpreted as a generalized sum of backward fluxes, where every local flux  has the same ``relevance''.

The following expression for the flux is a generalized sum of  weighted backward fluxes. There is now a kernel  which assigns  more weight (``importance'') to the nearest temperature gradients, that is:

\begin{equation}\label{Weighted-Average-FourierLaw}
J(x,t)=-\frac{\eta_\al}{\G(\al)}\int_{-\infty}^t (t-\tau)^{\al-1} k\frac{\p}{\p x}u(x,\tau) \dd\tau.
\end{equation}
Here, $\al$ is a constant in the interval $(0,1)$ that plays an important role, and $\eta_\al$ is a constant imposed to equate units of measures. Both will be specified later.  
      
Note that (\ref{FourierLaw}) and  (\ref{Weighted-Average-FourierLaw}) result from considering the kernels $K_1(t)\equiv \delta(t)$ and \linebreak $K_2(t)= \eta_{\al}\frac{t^{\al-1}}{\G(\al)}$, respectively, in the generalized flux equation $(\ref{Gen-FourierLaw})$.\\
	\bigskip
	
\item   \textbf{Implicit forms for the flux: $F(x,t,J(x,t))=G(x,t)$}.\\

One of the most famous formulations for the flux, is given by the  Cattaneo's equation \cite{Ca:1948}
\begin{equation}\label{Relaxed-FourierLaw}
J(x,t) + \tilde{\tau} \frac{\p}{ \p t} J(x,t)=-k\frac{\p}{\p x}u(x,t),
\end{equation}

 which was proposed with the aim of introducing an alternative to the ``unphysical'' property of the diffusion equation known as \textit{infinite speed of propagation}. Equation (\ref{Relaxed-FourierLaw})  can be seen as a first order Taylor approximation of (\ref{Taylor}) in which the flux is allowed to adjust to the gradient of the temperature according to a relaxation time $\tilde{\tau}$,
\begin{equation}\label{Taylor}
J(x,t+\tilde{\tau} )=-k\frac{\p}{\p x}u(x,t).
\end{equation}

Another approach assumes that the integral of the back fluxes, at the current time, is proportional to the gradient of the temperature:
\begin{equation*}\label{Imp-Average-J}
\frac{1}{\tilde{\tau}}\int_{-\infty}^t J(x,\tau)\dd\tau= -k\frac{\p}{\p x}u(x,\tau).
\end{equation*}

Yet another formulation considers that \textsl{the integral of the weighted backward fluxes at the current time, is proportional to the gradient of the temperature}:
\begin{equation}\label{int-alpha-J}
\frac{\nu_\al}{\G(1-\al)}\int_{-\infty}^t (t-\tau)^{-\al}J(x,\tau) \dd\tau=-k \frac{\p}{\p x}u(x,\tau).
\end{equation}

\begin{nota} Although when we talk about backward  fluxes  it is logical to consider the lower limit of the integral at $-\infty$, we can suppose that the function $u$ has   remained constant (for some reason) for all  $t <0$, where with $0$ we refer to a certain initial time.  Moreover, under this condition, that is $ u(x,t) \equiv u_0 $, for every $ t <0$, the expressions (\ref{Weighted-Average-FourierLaw}) and (\ref{int-alpha-J}) become
\begin{equation}\label{Weighted-Average-FourierLaw_0}
J(x,t)=-\frac{\eta_\al}{\G(\al)}\int_{0}^t (t-\tau)^{\al-1} k\frac{\p}{\p x}u(x,\tau) \dd\tau,
\end{equation}

and 
\begin{equation}\label{int-alpha-J_0}
\frac{\nu_\al}{\G(1-\al)}\int_{0}^t (t-\tau)^{-\al}J(x,\tau) \dd\tau=-k \frac{\p}{\p x}u(x,\tau),
\end{equation}

respectively.
\end{nota}

Expressions (\ref{Weighted-Average-FourierLaw_0}), (\ref{int-alpha-J_0}) are closely linked to fractional calculus. Let us present the basic definitions that will be employed throughout the article.
 
\end{itemize} 
 
\begin{defi}\label{defi frac} Let $\left[a,b\right]\subset \bbR$ and $\al \in \bbR^+$ be such that $n-1<\al\leq n$. 
\begin{enumerate}
	\item   For $f \in L^1[a,b]$, we define the \textsl{fractional 
Riemann--Liouville integral of order  $\alpha$} as
$$_{a}I^{\alpha}f(t)=\frac{1}{\Gamma(\alpha)}\int^{t}_{a}(t-\tau)^{\alpha-1} f(\tau)\dd\tau. $$
\item For $f\in AC^n[a,b]=\left\{ f \;|\;  f^{(n-1)} \, \text{ is absolutely continuous on [a,b]} \right\}$, we define the 
\textsl{fractional Riemann--Liouville  derivative of order $\alpha$ } as
$$ ^{RL}_{a}D^{\alpha}f(t)= \left[ D^n \, _{a}I^{n-\alpha}f  \right] (t) =\frac{1}{\Gamma(n - \alpha)}\frac{d^n}{dt^n}\int^{t}_{a}(t-\tau)^{n-\alpha-1} f(\tau)\dd\tau. $$
\item For $ f \in W^n(a,b)=\left\{ f  | \, f^{(n)}\in L^1[a,b] \right\}$, we define the \textsl{fractional Caputo derivative of order  $\al$} as
$$\,^C_{a} D^{\alpha}f(t)= \left[ \, _{a}I^{n-\alpha}(D^nf)  \right] (t) =\left\{\begin{array}{lc} \frac{1}{\Gamma(n-\al)}\displaystyle\int^{t}_{a}(t-\tau)^{n-\al-1} f^{(n)}(\tau)\dd\tau, &  n-1<\al<n,\\
f^{(n)}(t), &   \al=n. \end{array}\right.$$
\end{enumerate}
\end{defi}

\begin{nota} With these definitions, equations (\ref{Weighted-Average-FourierLaw_0}) and (\ref{int-alpha-J_0}) can be rewritten as

\begin{equation*}\label{Weighted-Average-FourierLaw-2}
J(x,t)=-\eta_\al k \;\,_0I^\al_t\frac{\p}{\p x}u(x,t) 
\end{equation*}
   and
\begin{equation*}\label{int-alpha-J-2}
\nu_\al \,_0I^{1-\al}_t J(x,\tau) =-k\frac{\p}{\p x}u(x,t).
\end{equation*}

\end{nota}

\noindent Table 1 exhibits the governing equations derived from (\ref{cont eq}) for the different choices of the flux $J$.
\begin {table}[H]\label{tablaFluxesEquations}
\caption {Flux and diffusion equations} 
\begin{flushleft}
\begin{tabular}{| l| l| l|} 
\hline& & \\
  \textbf{Equation for the Flux} & \textbf{Resulting Diffusion   Eq.} & \textbf{Observations}  \\ [1ex] 
 \hline  
 $J=-k\frac{\p u}{\p x}$  & $\frac{\p u}{\p t}= d\frac{\p^2  u}{\p x^2}$ & Heat equation   \\ [1ex] 
\hline  
 $J=-\frac{1}{\tilde{\tau}}\int_0^t k\frac{\p u}{\p x} d\tau$     &  $\frac{\p^2 u}{\p t^2} =\frac{d}{\tilde{\tau}} \frac{\p^2 u}{\p x^2}$  &  Wave equation  \\ [1ex] 
\hline 
 $J=\,-k \eta_\al\, _0I^\al_t \frac{\p u}{\p x}$  & $ \,^C_0D^{1+\al}_t u=\eta_\al d \frac{\p^2 u}{\p x^2}$ & Superdiffusion equation\\ [1ex] 
\hline 
 $J + \tilde{\tau} \frac{\p}{\p t} J=-k\frac{\p u}{\p x}$   &  $\frac{\p u}{\p t}+ \tilde{\tau} \frac{\p^2 u}{\p t^2}=  \frac{\p^2 u}{\p x^2}$ & Telegraph equation \\ [1ex] 
\hline 
 $\frac{1}{\tau^*}\int_0^t J \dd\tau= -k\frac{\p u}{\p x} $ & $\frac{\p }{\p t}\left(u-\frac{d}{\tau^*}\frac{\p^2 u}{\p x^2}\right) =0 $  & Elliptic equation \\ 
 & & with parameter $t$\\ 
\hline 
 $\nu_\al \,_0I^{1-\al}_t J =-k\frac{\p u}{\p x}$ & $^C_0D^\al u(x,t)=\frac{d}{\nu_\al}\frac{\p^2}{\p x^2}u(x,t)$ & Subdiffusion equation \\ [1ex]
\hline
\end{tabular}
\end{flushleft}
\end{table}
 There are  many references  about these different fluxes and their corresponding  governing equations \cite{Cannon,Ca:1948,Eidelman-Kochubei-LIBRO,FM-Libro,Povstenko,Pskhu-Libro}. Specially, the subdiffusion equation is one of the most studied in the last 10 years: The Cauchy problem \cite{GoLuMa:2000, EidelmanKochubei:2004, LuMaPa-TheFundamentalSolution, Povs:2008}, initial and boundary value problems \cite{GoReRoSa:2015, SaYa:2011}, maximum principles \cite{Luchko:2009, AlRefai-Luchko, Liu:2017, Ro:2016}. Nevertheless,  fractional phase change problems have been very poorly studied \cite{At:2012, VMD:2012}. Some of these articles propose a physical approach \cite{Ce:2018, GeVoMiDa:2013, Voller:2010, VoFaGa:2013, Voller:2014} and others do a purely mathematical treatment \cite{Li-Xu-Jiang, RoSa:2013, RoTa:2014}. \\
 The goal of this paper is to present a new mathematical model for a  one phase change problem with a memory flux, which derives  in a fractional free boundary problem, such that the governing  equations of this model are consistent both mathematically and physically speaking. We will pay special attention to the interchange of limits and integrals, which is a sensitive issue when working with fractional derivatives (see \cite{RoTa:2017-TwoDifferent}).\\
 In Section 2, some properties  of fractional calculus which will be useful later are provided.\\
 In Section 3, a mathematical formulation  for an instantaneous phase-change problem for a material with memory is presented. In this model, an implicit equation for the flux involving fractional integrals  is used.

 Finally, in Section 4, an equivalent formulation is presented, which allows us to give an integral relation for the free boundary, which we consider important in future research on existence and  uniqueness of solutions, or properties of the free boundary.

\section{Preliminaries of Fractional Calculus}

\begin{propo}\label{propo frac}\cite{Diethelm} The following properties involving the fractional integrals and derivatives hold:
\begin{enumerate}
\item \label{RL inv a izq de I} The  \textsl{fractional Riemann--Liouville derivative } is a left inverse operator of the \textsl{fractional Riemann--Liouville integral} of the same order  $\al\in \bbR^+$. If $f \in L^1[a,b]$, then
$$^{RL}_{a}D^{\al}\,_{a}I^{\al}f(t)=f(t)  \quad a.e.$$

\item The fractional  Riemann--Liouville integral is not, in general, a left inverse operator of the fractional derivative of Riemann--Liouville.\\
 
In particular, if $0<\al<1$, then 
$ _{a}I^{\al}(^{RL}_{a}D^{\al}f)(t)=f(t) - \dfrac{_{a}I^{1-\al}f(a^+)}{\G(\al)(t-a)^{1-\al}}.$ 

\item\label{caso part I inv de RL} If there exists some $\phi \, \in L^1(a,b)$ such that $f=\,_aI^\al \phi$, then 
$$_{a}I^{\al}\,^{RL}_{a}D^{\al} f(t)=f(t)  \quad a.e.$$

\item\label{relacion RL-C} If $n-1<\al\leq n$ and  $f\in AC^n[a,b],$ then   $$^{RL}_{a}D^{\al}f (t)=\displaystyle\sum_{k=0}^{n-1} \frac{f^{(k)}(a)}{\G(1+k-\al)}(t-a)^{k-\al}+\, ^C_{a} D^{\alpha}f(t).$$
In particular, for $0<\al<1$, we have 
$$ ^{RL}_{a}D^{\al}f (t)=\frac{f(a)}{\G(1-\al)}(t-a)^{-\al}+\, ^C_{a} D^{\alpha}f(t).$$

\end{enumerate}

\end{propo}

\begin{propo} \cite{Samko} The following limits hold:
\begin{enumerate}
\item If we set $_aI^0=Id$, the identity operator, then for every $ f \, \in L^1[a,b]$,   
$$ \displaystyle\lim_{\al\searrow 0}\, _aI^\al f(t) =_aI^0f(t)=f(t). $$

\item  For every $f \in C^1(a,b)$, 
$$ \displaystyle\lim_{\al\nearrow 1}\,_a^CD^\al f(t) = f'(t) \quad \text{and} \quad  \displaystyle\lim_{\al\searrow 1}\,_a^CD^\al f(t) = f'(t)-f'(0^+),\quad \forall t\in [a,b]. $$

\item  For every $f \in AC^1[a,b]$, 
$$ \displaystyle\lim_{\al\nearrow 1}\,_a^{RL}D^\al f(t) = f'(t)\quad\text{and} \quad \displaystyle\lim_{\al\searrow 1}\,_a^{RL}D^\al f(t) = f'(t),\quad a.e. \; t\in (a,b). $$

\end{enumerate}

\end{propo}

\begin{obs} If we consider a function $f$ supported in $[0,\infty)$ and $\chi_\alpha$ is the locally integrable function defined by 
\begin{equation*}\label{chi_alpha}\chi_\alpha(t)=\begin{cases} \frac{t^{\alpha-1}}{\G(\al)} & \text{if } t>0, \\
0 & \text{if } t\leq 0,\end{cases}
\end{equation*}
 then we have the following properties for $0<\al<1$:\\
\begin{equation*}\label{I_al con convolucion}
 _{a}I^{\al}f(t)=\left(\chi_\alpha \ast f\right) (t),
\end{equation*}
\begin{equation*}\label{D_al RL con convolucion}
 _{a}^{RL}D^{\al}f(t)=\frac{d}{dt}\left(\chi_{1-\alpha} \ast f\right) (t),
\end{equation*}
\begin{equation*}\label{D_al C con convolucion}
 _{a}^CD^{\al}f(t)=(\chi_{1-\alpha}) \ast \frac{d}{dt}f (t).
\end{equation*}

\end{obs}

\section{Modelling a Phase Change Problem with a flux with memory: A fractional Stefan problem}\label{FSP}

The aim of this section is to formulate  mathematical models associated to a  one--dimensional fractional phase change problem.

The classical phase change problems for the heat equation obtained by considering the Fourier Law for the flux are known in the mathematical literature as free boundary problems, and under certain conditions as Stefan problems.

The fundamental equations involved in Stefan problems are: the heat equation and the Stefan condition (derived from the connection between  the velocity of the free boundary and the heat fluxes of the two temperatures corresponding to the different phases).

We will focus on deriving the fractional diffusion equation  and (making an abuse of language) the ``fractional Stefan condition''.
\medskip
\medskip

\noindent\textbf{Physical problem:} Melting of a  semi--infinite slab  ($0\leq x<\infty$) of  a material with memory, which is at the melt temperature  $T_m$, by imposing a constant temperature $T_0>T_m$ on the fixed face $x=0$. All the thermophysical parameters are constants.  \\

\noindent \textbf{Mathematical problem}  Let $u=u(x,t)$ be the temperature and  let $J(x,t)$ be the memory flux of the material at position $x$ and time $t$. Let $x=s(t)$ be the function representing the (unknown) position of the free boundary at time $t$ such that $s(0)=0$. We will assume that $s$ is an increasing function and consequently, an invertible function.\\

The flux modelling the material with  memory is considered under the assumption  that \textsl{the generalized sum of the weighted backward fluxes at the current time is proportional to the gradient of the temperature}, that is
\begin{equation}\label{J-elegido}
\nu_\al\,_{h(x)}I^{1-\al}_tJ(x,t)=-k\frac{\p u}{\p x}(x,t),
\end{equation}
where the initial time in the fractional integral is given by the function  $h$ which gives us the time when the phase change occurs. That is
$$ h(x)=s^{-1}(x) \qquad (\text{i.e. } \,\, x=s(t) ).$$

\begin{center}
\begin{figure}[h]
\hspace{4cm}\includegraphics[width=0.5\textwidth]{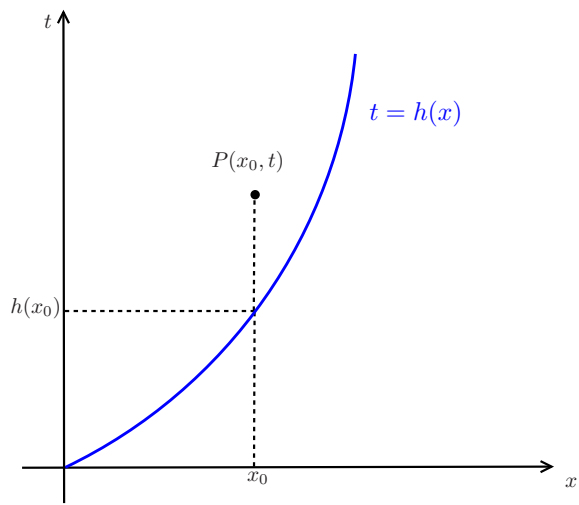}
\caption{The free boundary $h(x)$ vs $x$.}
\end{figure}
\end{center}

The parameter  $\nu_\al$ is a parameter with physical dimension such that 
\begin{equation}\label{lim-nu-al}
\displaystyle\lim_{\al\nearrow 1}\nu_{\al}=1.
\end{equation}
 
This parameter  has been added to preserve the consistency with respect to the units of measure in equation  (\ref{J-elegido}).  In fact, considering the units of measure given in (\ref{medidas}), we have
 \begin{equation}\label{med-J}
 \left[ J \right]= \left[ k u_x \right]=\frac{\textbf{m}}{\textbf{t}^3}, 
 \end{equation}
 
 \begin{equation}\label{med-IalJ}
 \left[ _{h(x)}I^{1-\al}_tJ(x,t) \right]=\left[\frac{1}{\G(1-\al)} \int_{h(x)}^t\frac{J(x,\tau)}{(t-\tau)^\al}\dd\tau \right]=\frac{\textbf{m}}{\textbf{t}^3}\frac{1}{\textbf{t}^\al}\textbf{t}=\frac{\textbf{m}}{\textbf{t}^{2+\al}}.  
 \end{equation}

Then, by  (\ref{med-J}), (\ref{med-IalJ}) and (\ref{J-elegido}) one gets

\begin{equation}\label{med-nu_al}
\left[\nu_\al\right]=\frac{\left[k\frac{\p u}{\p x}\right]}{\left[_{h(x)}I^{1-\al}_tJ\right])}=\frac{1}{\textbf{t}^{1-\al}}.
\end{equation}

\begin{obs}  Due to the properties of the Riemann--Liouville integral, the limit expression for $\al=1$  in  (\ref{J-elegido}) yields  the classical Fourier Law.  
\end{obs}

\begin{obs}\label{nota flujo nulo} Notice that, since we are assuming that the temperature is constant for $x>s(t)$, then the gradient of the temperature  is null in the region $x>s(t), t>0$, which implies that 
\begin{equation}\label{null-flux}
\nu_\al\,_{0}I^{1-\al}_tJ(x,t)=0, \qquad \forall \, x>s(t),\, t>0. 
\end{equation}
Applying the inverse operator $ ^{RL}_0D^{1-\al}_t $ to both sides of equation (\ref{null-flux}) leads to

\begin{equation*}\label{null-flux-2}
\nu_\al J(x,t)=\,^{RL}_0D^{1-\al}_t 0=0 , \qquad \forall \, x>s(t),\, t>0. 
\end{equation*}
Then, for every $(x,t)$ such that $0<x<s(t)$, $t>0$ it results that 
\begin{equation*}\label{I_0=I_h}
\begin{array}{rcl}
 \nu_\al\, _{0}I^{1-\al}_tJ(x,t) & = & \dfrac{\nu_\al}{\G(1-\al)} \displaystyle\int_{0}^t\frac{J(x,\tau)}{(t-\tau)^\al}\dd\tau \\[1ex]
  & =  & \dfrac{\nu_\al}{\G(1-\al)} \displaystyle\int_{0}^{h(x)}\frac{0}{(t-\tau)^\al}\dd\tau + \dfrac{\nu_\al}{\G(1-\al)} \displaystyle\int_{h(x)}^t\frac{J(x,\tau)}{(t-\tau)^\al}\dd\tau \\ [1ex]
  & =  & \nu_{\al}\,_{h(x)}I^{1-\al}_tJ(x,t).\end{array}\end{equation*}

\noindent So, assuming that $u(x,t)\equiv T_m$ in the region $x>s(t),\, t>0$, condition (\ref{J-elegido}) is equivalent to 

\begin{equation}\label{J-elegido-equiv}
\nu_\al\,_{0}I^{1-\al}_tJ(x,t)=-k\frac{\p u}{\p x}(x,t) \qquad \forall \, 0<x<s(t),\, t>0.
\end{equation}

 However,  in the following,  expression (\ref{J-elegido}) will be chosen   since  the dependence on starting time (linked to the free boundary) may be overlooked if we consider (\ref{J-elegido-equiv}).    
\end{obs}

Now, being the Riemann--Liouville fractional derivative of order $1-\al$ a left inverse operator of the fractional Riemann--Liouville integral (Proposition $\ref{propo frac}-\ref{RL inv a izq de I})$, an explicit expression for the memory flux at position $x$ and time $t$ can be derived, and it is given by

\begin{equation}\label{J-explicit}
J(x,t)=-\frac{k}{\nu_\al}\,_{h(x)}^{RL}D^{1-\al}_t \frac{\p u}{\p x}(x,t),
\end{equation} 
or
\begin{equation}\label{J-explicit-2'}
J(x,t)= -\frac{k}{ \nu_{\al}}\frac{1}{\G(\al)}\frac{\p}{\p t} \int_{h(x)}^t (t-\tau)^{\al-1}\frac{\p u}{\p x}(x,\tau)\dd\tau.
\end{equation}

Putting  
\begin{equation}\label{mu_al}
\mu_\al=\frac{1}{\nu_\al},
\end{equation}
 from (\ref{med-nu_al}) and (\ref{lim-nu-al}) it results that
\begin{equation*}\label{med-mu-al}
[\mu_\al]=\textbf{t}^{1-\al}
\end{equation*}
 and 
\begin{equation}\label{mu->1}
\displaystyle\lim_{\al\nearrow 1}\mu_{\al}=1.
\end{equation}
Then, equation (\ref{J-explicit-2'}) becomes
\begin{equation}\label{J-explicit-2}
J(x,t)= -k\mu_{\al}\,_{h(x)}^{RL}D^{1-\al}_t \frac{\p u}{\p x}(x,t).
\end{equation}

\begin{obs} Fractional explicit expressions for the flux such as  the given in  (\ref{J-explicit-2}) were considered in many publications (see for instance \cite{Povstenko, RoTa:2017}). Although it is a direct consequence of the formulation (\ref{J-elegido}), up to  now, the physical meaning of the partial derivative with respect to time in (\ref{J-explicit-2}) was not clear.
\end{obs}

Let us derive the governing equations of our problem. Note that the starting time being a function of $x$  in the fractional derivative, the governing equation will not coincide exactly with the subdiffusion equation given in Table \ref{tablaFluxesEquations}.\\

Let $0<x<s(t), t>0$ be. Differentiating equation (\ref{J-elegido}) respect to $x$ yields that
$$
\frac{\p}{\p x}\left(\, \nu_\al\,_{h(x)}I^{1-\al}_tJ(x,t)\right)=-k\frac{\p^2 u}{\p x^2}(x,t).
$$
Or equivalently, 
\begin{equation*}\label{gov-eq-1}
\frac{\p}{\p x}\left[ \frac{\nu_\al}{ \G(1-\al)} \int_{h(x)}^t (t-\tau)^{-\al}J(x,\tau)d\tau \right] =-k\frac{\p^2 u}{\p x^2}(x,t).\end{equation*}
Differentiating  the left-hand side of latter equation and using the continuity equation (\ref{cont eq}) we have 
\begin{equation}\label{gov-eq-2}
 \frac{\nu_\al}{ \G(1-\al)} \int_{h(x)}^t(t-\tau)^{-\al} \rho c \frac{\p}{\p t}u(x,\tau)d\tau +\nu_\al \displaystyle\lim_{\tau\searrow h(x)}\frac{(t-\tau)^{-\al}}{ \G(1-\al)} J(x,\tau)h'(x) =k\frac{\p^2 u}{\p x^2}(x,t).
 \end{equation}

Then the governing equation is
\begin{equation}\label{gov-eq-3}
 \rho c \, _{h(x)}^CD^{\al}_t  u(x,t)+\displaystyle\lim_{\tau\searrow h(x)}\frac{(t-\tau)^{-\al}}{\G(1-\al)}J(x,\tau)h'(x) =\frac{ k }{\nu_\al}\frac{\p^2 u}{\p x^2}(x,t).\end{equation}

\begin{obs}
In case we use the alternative flux definition (\ref{J-elegido-equiv}), which is equivalent to (\ref{J-elegido}), in the derivation steps of the governing equation, we get 
$$\frac{\p}{\p x}\left[ \frac{\nu_\al}{ \G(1-\al)} \int_{0}^t (t-\tau)^{-\al}J(x,\tau)d\tau \right] =-k\frac{\p^2 u}{\p x^2}(x,t).$$
It must be pointed out that  the flux $J$ is not differentiable at $\tau=h(x)\in (0,t)$. Then, we can not differentiate under  integral in the left-hand side of the latter equation. This fact is the main reason why we will not arrive to a single Caputo derivative over  $[0,t]$ for $u$ in the left-hand side of  equation (\ref{gov-eq-3}), as has already  been suggested in literature.

\end{obs}

Now, we turn to study  the moving interface. The interface is a  curve where a discontinuity of the flux occurs. So, the energy balance  between  the latent heat and the difference of fluxes is given by  the Rankine--Hugoniot conditions at the interface 

\begin{equation}\label{cond de Stefan gral unidim}
\llbracket \textbf{J}  \rrbracket^s_l   =-\rho l \dot{s}(t).
\end{equation}
Here,  the double brackets represents the difference between the limits of the fluxes  from the solid phase and the liquid phase. Recall that the explicit flux is given by (\ref{J-explicit-2}) in the liquid phase,  and  the temperature is constant in the solid phase (which implies that the flux is null in this region as we have seen in Remark \ref{nota flujo nulo}). Then condition (\ref{cond de Stefan gral unidim}) becomes

\begin{equation*}\label{FractionalStefanCond-0}
\lim\limits_{x\nearrow s(t)}J(x,t) =\rho l  s'(t),  
\end{equation*}
or equivalently (by using (\ref{J-explicit-2}))

\begin{equation}\label{FractionalStefanCond}
-k\mu_\al \lim\limits_{x\nearrow s(t)}\, _{h(x)}^{RL}D^{1-\al}_t \frac{\p u}{\p x}(x,t) =\rho l  s'(t).
\end{equation}

Making an abuse of language, we will call equation (\ref{FractionalStefanCond})  the ``fractional Stefan condition''.

Assuming the continuity of the flux in the liquid region, the following equality holds
\begin{equation}\label{Cont del flujo}
\lim\limits_{x\nearrow s(t)}J(x,t) =\lim\limits_{t\searrow h(x)}J(x,t). 
\end{equation}
Combining (\ref{gov-eq-3}) and (\ref{Cont del flujo}) we get the following governing equation for the liquid phase
\begin{equation}\label{gov-eq-4}
\rho c _{h(x)}^CD^{\al}_tu(x,t)+\frac{\rho l}{\G(1-\al)}\frac{ s'(h(x)) h'(x) }{(t-h(x))^{\al}}=\frac{k}{\nu_\al } \frac{\p^2 u}{\p x^2}(x,t).\end{equation}
 Being $h$ the inverse function of $s$, it results that
\begin{equation}\label{s'h'}
h'(x)=\frac{1}{s'(s^{-1}(t))}=\frac{1}{s'(h(x))}.
\end{equation}
Finally, using (\ref{s'h'}) in (\ref{gov-eq-4}) leads to
\begin{equation}\label{gov-eq-posta}
\rho c\, _{h(x)}^CD^{\al}_tu(x,t)+\rho l \frac{(t-h(x))^{-\al}}{\G(1-\al)}=\frac{ k}{\nu_\al} \frac{\p^2 u}{\p x^2}(x,t).\end{equation}


 If we consider  the Stefan number defined by
 \begin{equation}\label{Ste num} Ste=\frac{c(T_0-T_m)}{l},   \qquad ([Ste]=1)
 \end{equation}
 
and we use it in (\ref{gov-eq-posta}), we get

\begin{equation}\label{gov-eq-posta-posta}
 _{h(x)}^CD^{\al}_tu(x,t)+ \frac{\left(T_0-T_m\right)}{Ste} \frac{(t-h(x))^{-\al}}{\G(1-\al)}=\mu_\al d \frac{\p^2 u}{\p x^2}(x,t),
 \end{equation}

where $d$ is the diffusion coefficient defined in  (\ref{medidas}) and $\mu_\al$ was given in (\ref{mu_al}).

\begin{nota}It is easy to check that 
$$ \left[\, _{h}^CD^{\al}_tu+ \frac{\left(T_0-T_m\right)}{Ste} \frac{(t-h(x))^{-\al}}{\G(1-\al)} \right]= \left[ \mu_\al d \frac{\p^2 u}{\p x^2}  \right]=\frac{T}{t^\al}.$$

\end{nota}

\begin{nota} We would like to highlight the difference between the fractional Stefan condition obtained in (\ref{FractionalStefanCond}) and the fractional Stefan condition considered in \cite{RoSa:2013} which was given by 
$$ \rho l  _{0}^CD^{\al}s(t)= -k \frac{\p u}{\p x}(s(t),t),$$
and was derived by replacing the classical derivative by the Caputo derivative in the classical Stefan condition. 
\end{nota}

\noindent Finally, using equations  (\ref{FractionalStefanCond}) and (\ref{gov-eq-posta-posta}), and adding appropriate initial conditions,  the system representing the physical problem proposed at the beginning of the current section is given by

\begin{equation}{\label{St-desde-0}}
\begin{array}{lll}
     (i)  &   _{h(x)}^CD^{\al}_tu(x,t)+ \frac{\left(T_0-T_m\right)}{Ste} \frac{(t-h(x))^{-\al}}{\G(1-\al)}=\mu_\al d \frac{\p^2 u}{\p x^2}(x,t), &   0<x<s(t), \,  0<t<T,  \, \,\\
     (ii) & s(0)=0,  \\ 
       (iii)  &  u(0,t)=T_0,   &  0<t\leq T,  \\
         (iv) & u(s(t),t)=T_m, & 0<t\leq T, \\
      (v) & \rho l s'(t)=- \mu_\al k \displaystyle\lim_{x\nearrow s(t)} \,^{RL}_{h(x)}D^{1-\al}_t \frac{\p}{\p x} u (x,t), & 0<t\leq T,
                                             \end{array}
                                             \end{equation}
where $h(x)=s^{-1}(x)$ for every $x>0.$

\begin{defi}\label{Def sol St-1} A pair $\{u,s\} $ is a solution of  problem $(\ref{St-desde-0})$  if the following conditions are satisfied
\begin{enumerate}
    \item $u$ is continuous in the region $\mathcal{R_T}=\left\{(x,t)\colon 0\leq x \leq s(t), \, 0<t\leq T \right\}$ and at  the point $(0,0) $, $u $ verifies that 
       $$ 0\leq \underset{(x,t)\rightarrow (0,0)}{\liminf}u(x,t)\leq \underset{(x,t)\rightarrow (0,0)}{\limsup } u(x,t)<+\infty .$$
	\item  $u\in $ $C(\mathcal{R_T})\cap C^2_x(\mathcal{R_T})$, such that $u \in 
	W^1_t((h(x),T))$ 	where $W^1_t((h(x),T)):=\{f(x,\cdot)\colon f \in W^1(h(x),T)  \quad \text{for every fixed } x\in [0,s(T)] \}$.
 \item $s \in C^1(0,T)$. 
		\item There exists $\left.^{RL}_{h(x)}D^{1-\al}_t \frac{\p}{\p x} u (x,t)\right|_{(s(t)^-,t)}$ for all $t \in (0,T]$.
    \item $u$ and $s$ satisfy $(\ref{St-desde-0})$.
   	
		\end{enumerate}
\end{defi}

\section{Integral condition}

It is interesting to note that, from the definition (\ref{J-elegido}) for the flux and Proposition $\ref{propo frac}-\ref{caso part I inv de RL}$, it results that expression (\ref{J-elegido}) is equivalent to expression  (\ref{J-explicit}) for the flux. Then, if we replace (\ref{J-explicit}) in the continuity equation (\ref{cont eq}) we obtain   the following governing equation, which is a fractional diffusion equation for the Riemann--Liouville derivative:

\begin{equation*}\label{gov-eq-RL}
\frac{\p u}{\p t}(x,t)= \mu_\al d \frac{\p}{\p x} \left(\,_{h(x)}^{RL}D^{1-\al}_t \frac{\p}{\p x}u(x,t)\right), \qquad 0<x<s(t), \, t>0. 
\end{equation*}   

\begin{lema}\label{Lema Salto}The following jumping formulas hold:
\begin{enumerate}
\item\label{1} If $ w(x,\cdot)$ and $ w_x(x,\cdot) \in L^1(0,T)$ then\\
 \begin{equation}\label{salto I} _{h(x)}I_t^{1-\al}\left[\frac{\p}{\p x} w(x,t)\right] - \frac{\p}{\p x} \left[\, _{h(x)}I_t^{1-\al} w(x,t)\right]  = \displaystyle\lim_{\tau \searrow h(x)}\,w(x,\tau)\frac{(t-\tau)^{-\al}}{\G(1-\al)}h'(x). \end{equation} 
\item\label{2} If $w(x,\cdot) \in AC^1(0,T) $ and  $ _{h(x)}^{RL}D^{1-\al}_t \left[ \frac{\p}{\p x} w(x,t) \right]$ is a continuous function then \\
 \begin{equation*}\label{salto D}   \frac{\p}{\p x} \left[ _{h(x)}^{RL}D^{1-\al}_t\, w(x,t)\right] -\, _{h(x)}^{RL}D^{1-\al}_t \left[ \frac{\p}{\p x} w(x,t) \right] = -\frac{\p}{\p t}\left(\displaystyle\lim_{\tau \searrow h(x)}\,w(x,\tau)\frac{(t-\tau)^{\al-1}}{\G(1-\al)}h'(x)\right). \end{equation*}
\end{enumerate}
\end{lema}
\proof
1. Applying first the definition of fractional integral  and differentiating  with respect to $x$  we get
\begin{multline}\label{salto1} \frac{\p}{\p x} \left[\, _{h(x)}I_t^{1-\al} w(x,t)\right]  =
\frac{\p}{\p x} \left[\, \frac{1}{\G(1-\al)}\int_{h(x)}^t w(x,\tau)(t-\tau)^{-\al}\dd \tau\right]  \\
 = \frac{1}{\G(1-\al)}\left[\int_{h(x)}^t \frac{\p }{\p x}w(x,\tau)(t-\tau)^{-\al}\dd \tau - \displaystyle\lim_{\tau \searrow h(x)}\,w(x,\tau)(t-\tau)^{-\al}h'(x)\right].
\end{multline}  
Equation (\ref{salto I}) can be derived directly from (\ref{salto1}).\\

\noindent 2. Analogously, 
\begin{multline*}\label{salto2} \frac{\p}{\p x} \left[ _{h(x)}^{RL}D^{1-\al}_t\, w(x,t)\right] =
\frac{\p}{\p x} \frac{\p }{\p t}\left[\,_{h(x)}I_t^{\al} w(x,t)\right]  \\
=\frac{\p }{\p t} \frac{\p}{\p x} \left[\,_{h(x)}I_t^{\al} w(x,t)\right]= \frac{\p }{\p t} \frac{\p}{\p x} \left[\frac{1}{\G(\al)}\int_{h(x)}^t w(x,\tau)(t-\tau)^{\al-1}\dd \tau \right]\\
 = \, _{h(x)}^{RL}D^{1-\al}_t \left[ \frac{\p}{\p x} w(x,t) \right] - \frac{\p }{\p t} \left(\displaystyle\lim_{\tau \searrow h(x)}\,\frac{w(x,\tau)(t-\tau)^{\al-1}}{\G(\al)}h'(x)  \right).
\end{multline*}

\endproof

\begin{propo} Consider the following fractional Stefan problem
\begin{equation}{\label{St-RL}}
\begin{array}{lll}
     (i)  & \frac{\p u}{\p t}(x,t)= \mu_{\al} d \frac{\p}{\p x} \left(\,_{h(x)}^{RL}D^{1-\al}_t \frac{\p}{\p x}u(x,t)\right), &   0<x<s(t), \,  0<t<T,  \, \,\\
     (ii) & s(0)=0,   & \\ 
       (iii)  &  u(0,t)=T_0>T_m,   &  0<t\leq T, \\
         (iv) & u(s(t),t)=T_m, & 0<t\leq T, \\
      (v) & \rho l s'(t)=-\mu_\al k \displaystyle\lim_{x\nearrow s(t)} \,^{RL}_{h(x)}D^{1-\al}_t \frac{\p}{\p x} u (x,t), & 0<t\leq T. 
                                             \end{array}
                                             \end{equation}
where $h$ is the function defined by  $h(x)=s^{-1}(x) $. Then problems (\ref{St-desde-0}) and (\ref{St-RL}) are \mbox{equivalent}.

\end{propo}

\proof  Being equations $(ii)$ to $(v)$ the same in both problems we have to check only that equations $(\ref{St-desde-0}-i)$ and 
$(\ref{St-RL}-i)$ are equivalent.   \\

Applying $_{h(x)}I^{1-\al}_t$ to both sides of $(\ref{St-RL}-i)$ we get

\begin{equation}\label{eq equiv 1}
  _{h(x)}^CD^{\al}_t u(x,t) = \mu_\al d \,  _{h(x)}I_t^{1-\al}\left[ \frac{\p}{\p x} \left(\,_{h(x)}^{RL}D^{1-\al}_t \frac{\p}{\p x}u(x,t)\right)\right].
\end{equation}

Proposition \ref{propo frac}--1 implies that if we apply $_{h(x)}^{RL}D^{1-\al}_t$ to both sides of ($\ref{eq equiv 1}$) we recover equation  $({\ref{St-RL}-i})$. Therefore $({\ref{St-RL}-i})$  and $(\ref{eq equiv 1})$ are equivalent.

On one hand, taking $w(x,t)= \,_{h(x)}^{RL}D^{1-\al}_t \left(\frac{\p}{\p x}u(x,t)\right)$  in Lemma \ref{Lema Salto}-\ref{1} we get  
\begin{multline}\label{eq nn}
\frac{\p}{\p x} \left( _{h(x)}I_t^{1-\al} \,_{h(x)}^{RL}D^{1-\al}_t \left(\frac{\p}{\p x}u(x,t)\right)\right) = \frac{\p}{\p x} \left[\frac{1}{\G(1-\al)}\int_{h(x)}^t (t-\tau)^{-\al}  \,_{h(x)}^{RL}D^{1-\al}_t\left( \frac{\p}{\p x}u(x,\tau)\right)\dd\tau \right]  \\
=  _{h(x)}I_t^{1-\al}\left[\frac{\p}{\p x} \left(\,_{h(x)}^{RL}D^{1-\al}_t \left(\frac{\p}{\p x}u(x,t)\right)\right)\right]  - \displaystyle\lim_{\tau \searrow h(x)}\,_{h(x)}^{RL}D^{1-\al}_t \left( \frac{\p}{\p x}u(x,t)\right)\frac{(t-\tau)^{-\al}}{\G(1-\al)}h'(x).
 \end{multline}

On the other hand, from  (\ref{J-elegido}) and Proposition $\ref{propo frac}-\ref{caso part I inv de RL}$, it holds that
\begin{equation}\label{eq nn2}
 _{h(x)}I_t^{1-\al}\left(\,_{h(x)}^{RL}D^{1-\al}_t \left( \frac{\p}{\p x}u(x,t)\right)\right)= \frac{\p}{\p x}u(x,t).
\end{equation}
Then $(\ref{eq nn})$ together with  $(\ref{eq nn2})$ yield
\begin{multline*}\label{eq equiv 2}
 \mu_\al d\, _{h(x)}I_t^{1-\al}\left[\frac{\p}{\p x} \left(\,_{h(x)}^{RL}D^{1-\al}_t \left(\frac{\p}{\p x}u(x,t)\right)\right)\right]  = \\
  = \mu_\al d\, \frac{\p^2}{\p x^2}u(x,t)-\frac{ d\,}{\G(1-\al)}\frac{\rho l }{ k} \frac{s'(h(x))h'(x)}{(t-h(x))^\al}.
 \end{multline*}

So, we can rewrite equation (\ref{eq equiv 1}) as 
\begin{equation}\label{eq equiv 3}
  _{h(x)}^CD^{\al}_t u(x,t) =  \mu_\al d\, \frac{\p^2}{\p x^2}u(x,t)-\frac{ d\,}{\G(1-\al)}\frac{\rho l }{ k} \frac{s'(h(x))h'(x)}{(t-h(x))^\al}.
\end{equation}

Taking into account $(\ref{s'h'})$ and $(\ref{Ste num})$, we conclude that $(\ref{eq equiv 3})$ is equivalent to  $(\ref{St-desde-0}-i)$, and then the thesis holds.
 
\endproof


\begin{defi}\label{Def sol St} A pair $\{u,s\} $ is a solution of  problem $(\ref{St-RL})$  if the following conditions are satisfied

\begin{enumerate}
    \item $u$ is continuous in the region $\mathcal{R_T}=\left\{(x,t)\colon 0\leq x \leq s(t), \, 0<t\leq T \right\}$ and at  the point $(0,0) $, $u $ verifies that 
      
 $$ 0\leq \underset{(x,t)\rightarrow (0,0)}{\liminf}u(x,t)\leq \underset{(x,t)\rightarrow (0,0)}{\limsup } u(x,t)<+\infty .$$
	\item  $u\in $ $C(\mathcal{R_T}^\circ)\cap C^2_x(\mathcal{R_T}^\circ)$, such that $u_x \in 
	AC^1_t((h(x),T))$ 	where $AC^1_t((h(x),T)):=\{f(x,\cdot)\colon f \in AC^1(h(x),T)  \quad \text{for every fixed } x\in [0,s(T)] \}$.
 \item $s \in C^1(0,T)$.
		\item There exists $\left.^{RL}_{0}D^{1-\al}_t \frac{\p}{\p x} u (x,t)\right|_{(s(t),t)}$ for all $t \in (0,T]$.
    \item $u$ and $s$ satisfy $(\ref{St-RL})$.
   	
		\end{enumerate}
\end{defi}

\begin{lema}\label{Relac-Deriv-RL-y-C-para una sol} If the pair $\left\{ u, s\right\}$ is a solution to problem $(\ref{St-RL})$ and $\frac{\p}{\p x}\left[ \, _{h(x)}^{RL}D_t^{1-\al}u(x,t)\right]$ is a continuous function, then 
$$ \frac{\p}{\p x} \, _{h(x)}^{C}D_t^{1-\al}u(x,t)= \, _{h(x)}^{RL}D_t^{1-\al}\left(\frac{\p}{\p x}u(x,t)\right).  $$
\end{lema}
\proof

Since  $\frac{\p}{\p x}\left[ \, _{h(x)}^{RL}D_t^{1-\al}u(x,t)\right]$ is a continuous function, the partial derivatives commutes and 

\begin{equation*}
\begin{split} \frac{\p}{\p x}\left[ \, _{h(x)}^{RL}D_t^{1-\al}u(x,t)\right] & =\frac{\p}{\p t} \frac{1}{\G(\al)} \frac{\p}{\p x} \int_{h(x)}^tu(x,\tau)(t-\tau)^{\al-1}d\tau \\
  & = \frac{\p}{\p t} \frac{1}{\G(\al)} \left[ \int_{h(x)}^t \left(\frac{\p}{\p x}u(x,\tau)\right)(t-\tau)^{\al-1}d\tau - u(x,h(x))(t-h(x))^{\al-1}h'(x)\right]\\
	& = \, _{h(x)}^{RL}D_t^{1-\al}\left(\frac{\p}{\p x}u(x,t)\right)- \frac{\p}{\p t}\frac{T_m}{\G(\al) }(t-h(x))^{\al-1}h'(x)\\
	& = \, _{h(x)}^{RL}D_t^{1-\al}\left(\frac{\p}{\p x}u(x,t)\right)+ \frac{T_m (1-\al)}{\G(\al) }\frac{h'(x)}{(t-h(x))^{2-\al}}\\
	& = \, _{h(x)}^{RL}D_t^{1-\al}\left(\frac{\p}{\p x}u(x,t)\right)+ \frac{\p }{\p x}\left[\frac{T_m }{\G(\al) (t-h(x))^{1-\al}}\right]
		\end{split}
	\end{equation*}
Then
\begin{equation}\label{CI-2}
\frac{\p}{\p x}\left[ \, _{h(x)}^{RL}D_t^{1-\al}u(x,t)-\frac{T_m }{\G(\al) (t-h(x))^{1-\al}}\right]= \, _{h(x)}^{RL}D_t^{1-\al}\left(\frac{\p}{\p x}u(x,t)\right).
\end{equation}
Applying  Proposition $\ref{propo frac}-\ref{relacion RL-C}$ in   (\ref{CI-2}) leads to
	
\begin{equation*}\label{CI-3}
\frac{\p}{\p x} \, _{h(x)}^{C}D_t^{1-\al}u(x,t)= \, _{h(x)}^{RL}D_t^{1-\al}\left(\frac{\p}{\p x}u(x,t)\right).
\end{equation*}

\noindent This concludes the proof.

\endproof

\begin{teo}\label{Teo cond de Stefan integral}
Let  $\{u,s\}$ be a solution of problem  $(\ref{St-RL})$ with $u$ such that $_{h(x)}^{RL}D_t^{1-\al} u(x,t)$ and  $ _{h(x)}^{RL}D_t^{1-\al} \left(\frac{\p}{\p x}u(x,t)\right)$ are in $\mathcal{C}^{1}\left( \mathcal{R}_T-\lbrace (0,0) \rbrace\right)$.  Then  the following integral relation for the free boundary $s(t)$ and the function $u(x,t)$  holds for every  $ \;t<T$:
\begin{equation}\label{cond de Stefan integral}
 \left( \frac{l}{c}-T_m\right)s^2(t)=2\mu_\al d \frac{T_0-T_m}{\G(\al+1)}t^\al -2\int_0^{s(t)}x u(x,t)dx-2 \mu_\al d \int_0^t\,\left. _{h(x)}^{C}D_t^{1-\al}u(x,t)\right|_{(s(\tau),\tau)}d\tau. 
\end{equation}
\end{teo}

\proof 
Recall the Green identity: 
$$ \int_{\p \Omega}Pdt+Qdx =\iint_{\Omega} \left(\frac{\p}{\p t}Q-\frac{\p}{\p x}P \right)\, dA , $$
where $\Omega$ is an open simply connected region,  $\p \Omega$ is a positively oriented, piecewise smooth, simple closed curve, and the field $F=(P,Q)$ is defined by
 \begin{equation*}\label{P}
\begin{split} P(x,t) & = -\mu_\al d\, x\,_{h(x)}^{RL}D_t^{1-\al} \left(\frac{\p}{\p x}u(x,t)\right) + \mu_\al d  \,\, _{h(x)}^{RL}D_t^{1-\al}(u(x,t)-T_m)\\
 & =- \mu_\al d \, x\,_{h(x)}^{RL}D_t^{1-\al} \left(\frac{\p}{\p x}u(x,t)\right)+\mu_\al d\,\, _{h(x)}^{C}D_t^{1-\al}u(x,t)
\end{split}
\end{equation*}
  
\begin{equation*}\label{Q} Q(x,t)=-x\, u(x,t).\end{equation*}  
Consider the region
${\mathcal{R}_t}^\e=\left\{(x,\tau) \in \bbR^2 \, / \, \e<\tau<t, 0<x<s(\tau)\right\}$ for $\e>0$ sufficiently small. Note that in this region, $F$ is $\CC^1$.
Now, taking into account that $u$ verifies $(\ref{St-RL}-i)$ and using Lemma \ref{Relac-Deriv-RL-y-C-para una sol}  we get
\begin{align*}\label{CI-1}
& \frac{\p}{\p t} Q (x,t)- \frac{\p}{\p x}P(x,t) = \nonumber \\
&=-x \frac{\p} {\p t}u(x,t)+\mu_\al d  \,_{h(x)}^{RL}D_t^{1-\al} \left(\frac{\p }{\p x}u(x,t)\right) \nonumber\\
&\qquad + \mu_\al d\, x\frac{\p}{\p x}\left[\,_{h(x)}^{RL}D_t^{1-\al} \left(\frac{\p }{\p x}u(x,t)\right) \right] - \mu_\al d \, \frac{\p}{\p x}\left[ \, _{h(x)}^{C}D_t^{1-\al}u(x,t) \right] \nonumber \\
&=-x\left[\frac{\p} {\p t} u(x,t)-\mu_\al d \, \frac{\p}{\p x}\left(\,_{h(x)}^{RL}D_t^{1-\al} \left(\frac{\p }{\p x}u(x,t)\right) \right)\right]\nonumber\\
&\qquad + \mu_\al d\,  \,_{h(x)}^{RL}D_t^{1-\al} \left(\frac{\p }{\p x}u(x,t)\right) - \mu_\al d \, \frac{\p}{\p x}\left[ \, _{h(x)}^{C}D_t^{1-\al}u(x,t) \right] \nonumber \\
&=0 \qquad \text{ for all } \, (x, t) \, \in \, \mathcal{R}_\e. 
\end{align*} 
  
Then, by Green's theorem one obtains
\begin{equation*}\label{int de linea cero}
 \int_{\partial \mathcal{R}_\e}Pd\tau+Qdx = 0.
\end{equation*}

Let  $
\partial \mathcal{R}_\e= \partial \mathcal{R}_1 \cup \partial \mathcal{R}_2 \cup \partial \mathcal{R}_3 \cup \partial \mathcal{R}_4
$ be, where $\partial \mathcal{R}_1  =  \left\{(x,\e) \colon 0\leq x\leq s(\e)\right\}$,\linebreak $\partial \mathcal{R}_2   = \left\{(s(\tau),\tau) \colon \e < \tau <t\right\}$,
$-\partial \mathcal{R}_3   =  \left\{(x,t) \colon 0\leq x\leq s(t)\right\}$ and $ -\partial \mathcal{R}_4  = \left\{(0,\tau) \colon \e \leq \tau \leq t\right\}.$\\

 Integrating  the field $(P,Q)$ over $\partial \mathcal{R}_\e$ we get

\begin{align}\label{IntGren}
&-\int_0^{s(\e)}x\,u(x,\e)\dd x -\int_\e^{t}\, \hspace{-0.2cm}\left[ \mu_\al d \, \left. _{h(x)}^{RL}D_t^{1-\al}\left( \frac{\p}{\p x}u(x,t)\right)\right|_{(s(\tau), \tau)}\right.\nonumber \\
&\left.+ \mu_\al d \left. _{h(x)}^{C}D_t^{1-\al}u(x,t)\right|_{(s(\tau), \tau)} \right] \dd \tau - \nonumber\\
&- \int_\e^{t} \,s(\tau) T_m s'(\tau)\dd \tau + \int_0^{s(t)}\,xu(x,t)\dd x - \int_\e^t \,\mu_\al d  _{h(x)}^{C}D_t^{1-\al} T_0 \,\dd \tau =0.
\end{align}
Applying the fractional Stefan condition (\ref{St-RL}) yields
\begin{equation}\label{CI-5}
\begin{split} 
-\int_0^{s(\e)}xu(x,\e)dx+ \left(\frac{l}{c}-T_m\right)\left[\frac{s(t)^2}{2}-\frac{s(\e)^2}{2}\right] + \mu_\al d \int_\e^t \left.  _{h(x)}^{C}D_t^{1-\al}u(x,t)\right|_{(s(\tau),\tau)}\\
+\int_{0}^{s(t)}xu(x,t)dx-\mu_\al d \frac{\left(T_0-T_m\right)}{\G(\al+1)}(t^{\al}-\e^\al)=0.
\end{split}
\end{equation}

Taking the limit when $\e \searrow 0$ in (\ref{CI-5}) it results that the integral relation (\ref{cond de Stefan integral}) holds as we wanted to prove.

\endproof

\begin{obs} It is worth noting the difference between
\begin{equation}\label{ojo-1}
\left.   _{h(x)}^{C}D_t^{1-\al}u(x,t)\right|_{(s(\tau),\tau)},
 \end{equation}
and
\begin{equation}\label{ojo-2} _{h(x)}^{C}D_t^{1-\al}u(s(t),t)= _{h(x)}^{C}D_t^{1-\al} T_m= 0. \end{equation}
In (\ref{ojo-1}) we first apply Caputo derivative and then evaluate at $(s(t), t)$. Instead, in (\ref{ojo-2}) we first evaluate function $u$ in $(s(t), t)$ and then the Caputo derivative is taken. 
 \end{obs}

\begin{obs} If we take $\al=1$, $T_m=0$ and all the physical constants equal to 1 in the integral relation (\ref{cond de Stefan integral}) we get
\begin{equation*}\label{cond de Stefan integral clasica}
 s^2(t)=-2\int_0^{s(t)}x u(x,t)dx + 2T_0\,t,
\end{equation*}
which is the classical integral relation for the free boundary when the classical Stefan problem is considered (see \cite{Cannon}--Lemma 17.1.1). \\
It was also   proved in \cite{Cannon} that  (\ref{cond de Stefan integral clasica}) is equivalent to the Stefan condition
\begin{equation}\label{cond de Stefan clasica}
s'(t)=-\frac{\p }{\p x}u(s(t), t), \quad \forall \, t>0. 
\end{equation} 
 Hence, it is natural to wonder if the ``fractional Stefan condition'' (\ref{FractionalStefanCond}) and the ``fractional integral relation'' (\ref{cond de Stefan integral}) are  equivalent as well.
\end{obs}

\begin{teo}\label{reciproco - cond de Stefan integral}
Let  $\{u,s\}$ be a solution of problem   $\left\{ (\ref{St-RL}-i), (\ref{St-RL}-ii), (\ref{St-RL}-iii), (\ref{St-RL}-iv), (\ref{cond de Stefan integral})\right\}$ such that $\frac{\p^2}{\p t\p x} u (x,t) \in \CC^1(\mathcal{R}_T) $,     $^{C}_{h(x)}D_t^{1-\al} u(x,t)|_{(s(t),t)} \in L^1(0,T)$. Then the functions  $s=s(t)$ and  $u=u(x,t)$ verify the fractional Stefan condition (\ref{FractionalStefanCond}).
\end{teo}
\proof

Reasoning as in Theorem \ref{Teo cond de Stefan integral}, we can state that 
again (\ref{IntGren}) holds. 

Taking the limit when $\e \searrow 0$ and using the integral relation (\ref{cond de Stefan integral}) it holds that 

\begin{equation}\label{teo 2-2}
 \frac{l}{c}s^2(t)= -2\mu_{\al} d \displaystyle\int_0^t\,s(\tau) ^{RL}D_t^{1-\al} \left. \frac{\p}{\p x}u (x,t)\right|_{(s(\tau),\tau)} d\tau .
\end{equation}
 Differentiating both sides of equation (\ref{teo 2-2}) whith respect to the $t-$variable and being \mbox{$s(t)>0$} for all $t>0$, the thesis holds.
\endproof

\section{Conclusions}

We have presented a physical phase change problem involving a material with memory. In the mathematical model,  a fractional Riemann--Liouville integral is used for an implicit definition of the flux. Then, the governing equations were obtained. As a result of this analysis two equivalent fractional Stefan problems (\ref{St-desde-0}) and (\ref{St-RL}) involving the Caputo and the  Riemann--Liouville derivative, respectively, were formulated. The comparison with the classical Stefan problem was given in each case. Moreover, the classical Stefan problem was recovered by making $\al \nearrow 1$. Finally, an integral relation which  is equivalent to the fractional Stefan condition was obtained.   

\section{Acknowledgements}

\noindent 

 The present work has been sponsored by the Projects PIP N$^\circ$ 0275 from CONICET--Univ. Austral, and ANPCyT PICTO Austral N$^\circ 0090$ (Rosario, Argentina). 
We appreciate the valuable suggestions by the anonymous referees which
helped  us to improve the paper. We are particularly grateful for the assistance given by Prof. Mar\'ia Soledad Aronna.

\bibliographystyle{plain}

\bibliography{Roscani_BIBLIO_GENERAL_nombres_largos2018_08}

\end{document}